\documentstyle[11pt]{article}
\onecolumn
\newtheorem{lem}{{\sc Lemma} }[section]
\newtheorem{thm}{{\sc Theorem} }[section]
\newtheorem{cor}{{\sc Corollary} }[section]

\newcommand{\e}{\equiv}
\newcommand{\nqv}{\equiv \!\!\!\!\!\!\//~~}

\title{\LARGE\bf Powersums representing residues mod $p^k$ ,\\
                        from Fermat to Waring }
\author{ {\sc Nico F. Benschop} \\[1ex]
  {\it benschop@chello.nl - Amspade Research} - The Netherlands, ~~11 March 2001 }
\date{}

\begin{document}
\maketitle

\begin{abstract}
The ring $Z_k(+,.)$ mod $p^k$ with prime power modulus (prime $p>2$) is
analysed. Its cyclic group $G_k$ of units has order $(p-1)p^{k-1}$, and
all $p$-th power $n^p$ residues form a subgroup $F_k$ with $|F_k|=|G_k|/p$.
The subgroup of order $p-1$, the $core~~A_k$ of $G_k$, extends Fermat's
Small Theorem ($FST$) to mod $p^{k>1}$, consisting of $p-1$ residues with
$n^p \e n$ mod $p^k$. The concept of {\bf carry}, e.g. $n'$ in $FST$ extension
$n^{p-1} \e n'p+1$ mod $p^2$, is crucial in expanding residue arithmetic to
integers, and to allow analysis of divisors of 0 mod $p^k$.
\\
For large enough $k \geq K_p$ (critical precison $K_p<p$ depends on $p$), all
nonzero pairsums of core residues are shown to be distinct, upto commutation.
The known $FLT$ case$_1$ is related to this, and the set $F_k+F_k$ mod $p^k$
of $p$-th power pairsums is shown to cover half of $G_k$. Yielding main result:
each residue mod $p^k$ is the sum of at most four $p$-th power residues.
Moreover, some results on the generative power (mod $p^{k>2}$) of divisors
of $p \pm 1$ are derived. ~~~[Publ.: {\it Computers and Mathematics with
Applications} V39 N7-8 (Apr.2000) p253-261]
\end{abstract}

{\bf MSC classes:} ~11P05, 11D41


\section{Introduction}

The concept of $closure$ corresponds to a mathematical operation composing two
objects into an object of the same kind. Structure analysis is facilitated by
knowing a minimal set of $generators$, to find preserved partitions viz.
congruences, that allow factoring the closure. For instance a finite state
machine decomposition using preserved (state-) partitions, corresponding to
congruences of the sequential closure (semigroup) of its state transformations.
\\ A minimal set of $generators$ is characterized by $anti ~closure$.
Then each composition of two generators produces a non-generator, thus
a new element of the closure. These concepts can fruitfully be used for
structure analysis of finite residue arithmetic.

  For instance positive integer $p$-th powers are closed under multiplication,
but no sum $a^p+b^p$ yields a $p$-th power for $p>$2 ~(Fermat's Last Theorem,
$FLT$).  Apparently $p$-th powers form an efficient set of additive generators.
Waring (1770) drew attention to the now familiar representation problem: the
sum of how many $p$-th powers suffice to cover all positive integers. Lagrange
(1772) and Euler showed that four squares suffice [2]. The general problem is
as yet unsolved.

Our aim is to show that four $p$-th power residues mod $p^k$ (prime $p>$2,
$k>0$ large enough) suffice to cover all $p^k$ residues under addition. As
shown in [3,4] the analysis of residues $a^p+b^p$ mod $p^k$ is useful here,
because under modulus $p^k$ the $p$-th power residues coprime to $p$ form a
proper multiplicative subgroup $F_k \e \{n^p\}$ mod $p^k$ of the group of units
$G_k(.)$ mod $p^k$, with $|F_k|= |G_k|/p$. The value range $F_k+F_k$ mod $p^k$
is studied.

Units group $G_k$, consisting of all residues coprime to $p$, is in fact known
to be cyclic for all $k>$0 [1]. There are $p^{k-1}$ multiples of $p$ mod $p^k$,
so its order $p^k-p^{k-1} =(p-1)p^{k-1}$ is a product of two coprime factors,
hence we have:

(1) ~~$G_k \e A_kB_k$ is a direct product of subgroups, with $|A_k|=p-1$
and $|B_k|=p^{k-1}$.

The $extension$ subgroup $B_k$ consists of all $p^{k-1}$ residues 1 mod $p$.
And in $core$ subgroup $A_k$, of order $|A_k|=p-1$ independent of $k$, each $n$
satisfies $n^p \e n$ mod $p^k$, denoted as $n^p \e n$. Hence core $A_k$ is
the extension of Fermat's Small Theorem ($FST$) mod $p$ to mod $p^k$ for $k>$1.
For more details see [4].

By a coset argument the nonzero corepairsums in $A_k+A_k$, for large enough $k$,
are shown to be all distinct in $G_k$, apart from commutation (thm2.1). This
leads to set $F_k+F_k$ of $p$-th power pairsums covering almost half of $G_k$,
the maximum possible in a commutative closure, and clearly related to Fermat's
Last Theorem ($FLT$) about the anti-closure of the sum of two $p$-th powers.

Additive analysis of the roots of 0 mod $p^2$, as sums of three $p$-th power
residues, via the generative power of divisors of $p \pm 1$ (thm3.1), yields
our main result (thm3.2): the sum of at most four $p$-th power residues
mod $p^k$ covers all residues, a {\it Waring-for-residues} result.
Finite semigroup- and ring- analysis beyond groups and fields is essential,
due the crucial role of divisors of zero.

\section{ Core increments as coset generators} 

The two component groups of $G_k \e A_k.B_k$ are residues mod $p^k$ of
two monomials: the $core$ function $A_k(n)=n^{q_k} ~(q_k=|B_k|=p^{k-1})$
and $extension$ function $B_k(n)=n^{|A_k|}=n^{p-1}$. Core function $A(n)$
has odd degree with a $q$-fold zero at $n$=0, and is monotone increasing
for all $n$. Its first difference $d_k(n)=A_k(n+1)-A_k(n)$ of even degree has
a global minimum integer value of 1 at $n=0$ and $n=-1$, and symmetry centered
at $n=-1/2$. Thus integer equality $d_k(m)=d_k(n)$ for $m \neq n$ holds only
if $m+n=-1$, called 1-complements.

Hence the next definition of a {\it critical precision $k=K_p$} for residues
with the same symmetric property is relevant for every odd $p$, not
necessarily prime. Core difference $d_k(n)$ is 1 mod $p$, so it is referred
to as {\bf core increment} $d_k(n)$. To simplify notation, the precision
index $k$ is sometimes omitted, with $\e$ denoting equivalence mod $p^k$,
especially since core $A_k$ has order $p-1$ independent of $k$.

{\bf Define} {\it critical precision} $K_p$ as the smallest $k$ for which
the {\it only} equivalences among the core-increments $d_k(n)$ mod $p^k$ are
the above described 1-complement symmetry for $n$ mod $p$, so these
increments are all distinct for $n=1~ .. ~(p-1)/2$.

Notice that $K_p$ depends on $p$, for instance $K_p$=2 for $p \leq 7$,
 ~$K_{11}=3, ~K_{13}=2$, and the next $K_p=4$ for $p=73$. Upperbound $K_p<p$
 will be derived in the next section (lem3.1c), so no 'Hensel lift' [7] occurs.
Notice that $|F_k|/|A_k|=p^{k-2}$, sothat $A_2=F_2=\{n^p\}$ mod $p^2$.\\

\begin{lem} 
Integer core-function $A_k(n)=n^{p^{k-1}}$ and its increment
      $d_k(n)=A_k(n+1)-A_k(n)$\\ \hspace*{1cm}
 both have period $p$ for residues mod $p^k$, with: \\
(a)~~ Odd symmetry $A_k(m) \e -A_k(n)$ at complements $m+n \e 0$ mod $p$, \\
(b)~~ Even symmetry $d_k(m) \e d_k(n)$ at 1-complements $m+n \e -1$ mod $p$.\\
(c)~~ Let $D_2$ be the set of distinct increments $d_2(n)$ mod $p^2$
 of $F_2=A_2$ for $0<n \leq (p-1)/2$, \\ \hspace*{1cm}
 then there are $|F_k+F_k ~\backslash ~0~|=|F_k|.|D_2|=|G_k|.|D_2|/p$ \\
\hspace*{1.5cm}  non-zero $p$-th power pairsums mod $p^k$ ~(any~ $k>$1).
\end{lem}
\begin{proof}
(a)~
Core function $A_k(n) \e n^{q_k}$ mod $p^k ~(q_k=p^{k-1}, ~n \nqv$ 0, -1 mod $p$)
has $p-1$ distinct residues for each $k>$0, satisfying $(n^q)^p \e n^q$ mod $p^k$,
with $A_k(n) \e n$ mod $p$ due to $FST$. Apparently, including $A_k(0) \e 0$ we have:
~$A_k(n+p) \e A_k(n)$ mod $p^k$ for each $k>$1, with {\bf period} $p$ in $n$. And
$A_k(n)$ of odd degree $q=q_k$ has {\bf odd symmetry} because:
\\ \hspace*{1cm} $A_k(-n) \e (-n)^q \e -n^q \e -A_k(n)$ mod $p^k$.

(b)~ Increment $d_k(n) \e  A_k(n+1)-A_k(n)$ mod $p^k$ also has period $p$ because
\\ \hspace*{1cm}
  $d_k(n+p) \e (n+p+1)^{q_k}-(n+p)^{q_k} \e (n+1)^{q_k}-n^{q_k} \e d_k(n)$ mod $p^k$.
\\
This yields residues 1 mod $p$ in extension group $B_k$. It is an even
degree polynomial, with leading term $q_k.n^{q_k -1}$, and {\bf even symmetry}:
\\ \hspace*{1cm}
   $d_k(n-1)=n^{q_k}-(n-1)^{q_k}=-(-n)^{q_k}+(-n+1)^{q_k}= d_k(-n)$,\\
   so $d_k(m)=d_k(n)$ mod $p^k$ for 1-complements: ~$m+n=-1$ mod $p$.

(c)~
Write $F$ for $F_k ~(any ~k>1)$, the subgroup of $p$-th power residues mod $p^k$
in units group $G_k$. Then subgroup closure $F.F=F$ implies $F+F=F.(F+F)=F.(F-F)$,
since $F+F=F-F$ due to $-1$ in $F$ for odd prime $p>$2. So non-zero pairsum set
$F+F ~\backslash ~0$ is the disjoint union of cosets of $F$ in $G$, as
generated by differences $F-F$. Due to (1):~ $G_k=A_k.B_k=F_k.B_k$, where
$A_k \subseteq F_k$, it suffices to consider only differences 1 mod $p$,
hence in extension group $B=B_k$, that is: in $(F-F) \cap B$. \\ This amounts
to $|D_2| \leq h=(p-1)/2$ distinct increments $d_2(n)$, for $n=1...h$
due to even symmetry (b), and excluding $n$=0 involving non-core $A_2(0)$=0.
These $|D_2|$ cosets of $F_k$ in $G_k$ yield:~ $|F_k+F_k ~\backslash~0~|=
|F_k|.|D_2|$, where $|F_k|=|G_k|/p=(p-1).p^{k-2}$ and $|D_2| \leq (p-1)/2$.
\end{proof}

For many primes $K_p=2$ so $|D_2|=(p-1)/2$, and Fermat's $p$-th power residue
pairsums cover almost half the units group $G_k$, for any precision $k>$1.
But even if $K_p>2$, with $|D_2|<(p-1)/2$, this suffices to express each
residue mod $p^k$ as the sum of at most four $p$-th power residues (thm3.2),
as shown in the next section.

\begin{thm} 
~For $a, b$ in core $A$ mod $p^k$, and $k \geq K_p$:  \\[1ex] \hspace*{1cm}
~All nonzero pairsums $a+b$ mod $p^k$ are distinct, apart from commutation,
so: \\[1ex] \hspace*{1in}
    ~$|~(A+A)~\backslash~0~|=~\frac{1}{2}~|A|^2=(p-1)^2 /2$.
\end{thm}
\begin{proof}
~Core $A_k$ mod $p^k ~(any ~k>1)$, here denoted by $A$ as subgroup of units
group $G$, satisfies $A.A \e A$ so the set of all core pairsums can be factored
as $A+A \e A.(A+A)$. Hence the nonzero pairsums are a (disjoint) union of the
cosets of $A$ generated by $A+A$.
Since $G \e A.B$ with $B=~\{n$=1 mod $p$\}, there are $|B|=p^{k-1}$ cosets of
$A$ in $G$. Then intersection $D \e ~(A+A) \cap B$ of all residues 1 mod $p$ in
$A+A$ generates $|D|$ distinct cosets of $A$ in $G$.

Due to $-1$ in core $A$ we have $A \e -A$ sothat $A+A \e A-A$. View set $A$ as
function values $A(n) \e n^{|B|}$ mod $p^k$, with $A(n) \e n$ mod $p ~(0<n<p)$.
Then successive core increments $d(n)=A(n+1)-A(n)$ form precisely
intersection $D$, yielding all residues 1 mod $p$ in $A+A \e A-A$. Distinct
residues $d(n)$ generate distinct cosets, so by definition of $K_p$ there
are for $k \geq K_p:~ |D|=(p-1)/2$ cosets of core $A$ generated by $d(n)$
mod $p^k$.
\end{proof}

\section{Core extensions from $A_k$ to $F_k$, and their pairsums mod $p^k$}

Extension group $B$ mod $p^k$, with $|B|=p^{k-1}$ has only subgroups of
order $p^e ~(e=0..k$-1). So $G \e A.B$ ~(1) has $k$ subgroups $X^{(e)}$
that contain core $A$, called {\bf core extensions}, of order
$|X^{(e)}|=(p-1).p^e$, with core $A=X^{(0)}, ~F=X^{(k-2)}$ and $G=X^{(k-1)}$.\\
Now $p+1$ generates $B$ of order $p^{k-1}$ in $G_k$ ~[4, lem2], and similarly:

(2) \hspace{1cm} $p^i+1$ of period $p^{k-i} ~~(i=1..k$-1) in $G$ generate
the $k-1$ subgroups of $B$.

Let $Y^{(e)} \subseteq B$, of order $p^e$, then all core extensions are cyclic
with product structure:

~~~~~~~~~~~~$X^{(e)} \e A.Y^{(e)}$ ~in $G(.)$
    ~where $|A|$ and $|Y^{(e)}|$ are relative prime.

Using (2) with $k-i=e$ yields: \\[1ex]
(2') \hspace{1.5cm}
 ~~~$Y^{(e)} \e (p^{k-e}+1)^* \e \{m.p^{k-e}+1\}$ mod $p^k$ ~~(all $m$).

As before, using residues mod $p^k$ for any $k>1:~ D \e (A-A) \cap B$ contains
the set of core increments. Then thm2.1 on core pairsums $A+A$ is generalized
as follows (lem3.1a) to the set $X+X$ of core extension pairsums
mod $p^j~(j>$1), with $F+F ~(Fermat~sums)$ for $j=k-2$.

Extend Fermat's Small Theorem $FST:~ n^{p-1} \e 1$ mod $p$ ~to~ $n^{p-1} \e n'p+1$
mod $p^2$, which defines the {\bf $FST$-carry} $n'$ of $n<p$. This yields  an
efficient {\bf core generation} method (b) to compute $n^{p^i}$ mod $p^{i+1}$,
as well as a proof (c) of critical precision upperbound $K_p <p$.\\

\begin{lem} .\\ \hspace*{.5cm}
For core increments $D_k=(A_k-A_k) \cap B_k$ in $G_k=A_k.B_k$ mod $p^{k>1}$
~(prime $p>2$), \\ \hspace*{.5cm}
     $p$-th power residues set $F_k \e \{n^p\}$ mod $p^k$,
      and $X_k$ any core extension $A_k \subseteq X_k \subseteq F_k$:\\
(a)~~~
 $X_k+X_k \e X_k.D_k$, so core-increments $D_k$ generate the $X_k$-cosets
 in $X_k+X_k$.\\
(b)~~~ $[n^{p-1}]^{p^{i-1}} \e n'.p^i+1$ ~~~mod $p^{i+1}$,
~~~where $FST$-carry $n'$ of $n$ does not depend on $i$,\\ \hspace*{1cm}
and:~ $n^{p^i} \e [n'p^i+1]n^{p^{i-1}}$ ~mod $p^{i+1}$.
\\
(c)~~~ For $k=p$~:~~ $|D_p|=(p-1)/2$ mod $p^p$, so critical precision $K_p<p$.
\end{lem}
\begin{proof}
(a)~
Write $X$ for $X^{(e)}_k$ then, as in theorem 1.1:  ~$X+X=X-X=(X-X)X$.
For residues mod $p^k$ we seek intersection $(X-X) \cap B$ of all distinct
residues 1 mod $p$ in $B$ that generate the cosets of $X$ in $X+X$ mod $p^k$.
By (2,~2') core extension $X= A.Y= A.\{mp^{k-e}+1\}$. Discard terms divisible
by $p$ (are not in $B$) then: $(X+X) \cap B = (A+A) \cap B=(A-A) \cap B=D$
for each core extension. So $A+A$ and $X+X$ have the same coset generators
in $G_k$, namely the core increment set $D=D_k \subset B_k$.

(b)~
Notice successive cores satisfy by definition $A_{i+1} \e A_i$ mod $p^i$.
In other words, each $p$-th power step $i \rightarrow i+1:~ [n^{p^i}]^p$
produces one more significant digit ($msd$) while fixing the $i$ less
significant digits ($lsd$). Now $n^{p-1} \e n'p+1$ mod $p^2$ has $p$-th
power residue ~$[n^{p-1}]^p \e n'p^2+1$ mod $p^3$, implying lemma part (b)
by induction on $i$ in $[n^{p-1}]^{p^i}$.\\
This yields an efficient {\bf core generation} method.
Denote $f_i(n) \e n^{p^i}$, with $n<p$, then:\\[1ex] \hspace*{1cm}
(3)~~~~~ $f_i(n) \e n^{p^i} \e [n^p]~^{p^{i-1}} \e [n.n^{p-1}]~^{p^{i-1}} \e
     f_{i-1}(n).[n'p^i+1]$ ~mod $p^{i+1}$, ~implying:\\ \hspace*{1cm}
(3')~~~~ $f_i(n) \e f_{i-1}(n)$ mod $p^i$,
    ~next core msd~ $f_{i-1}(n)n'p^i \e nn'p^i \nqv 0$ mod $p^{i+1}$.

Notice that by $FST:~~ f_k(n) \e n$ mod $p$ for all $k \geq 0$,
and $0<n<p$ implies $n' \nqv 0$ mod $p$.

(c)~
In (a) take $X_k=F_p$ and notice that $F_p+F_p \e F_p-F_p$ mod $p^p$ contains
$h$ distinct integer increments $e_1(n)=(n+1)^p-n^p < p^p$ ..(4).. which are
1 mod $p^p$, hence in $B_p$. They generate $h$ distinct cosets of core $A_p$
in $G_p \e A_pB_p$ mod $p^p$, although they are not core $A_p$ increments.
Repeated $p$-th powers $n^{p^i}$ in constant $p$-digit precision yield
increments $e_i(n) \e (n+1)^{p^i}-n^{p^i}$ mod $p^p$, which for $i=p-1$ produce
the increments of core $A_p$ mod $p^p$. \\ Distinct increments $e_i(n) \nqv
e_i(m)$ mod $p^p$ remain distinct for $i \rightarrow i+1$, shown as follows.

For non-symmetric $n,~m<p$ (lem2.1b) let increments $e_i$ satisfy:

~~~~~(5)~~~~~ $e_i(n) \e e_i(m)$ ~mod $p^j$ ~for some $j<p$, ~and~

~~~~~(5')~~~~ $e_i(n) \nqv e_i(m)$ ~mod $p^{j+1}$.

Then for $i \rightarrow i+1$ the same holds, since $e_{i+1}(x)=
[f_i(x+1)]^p-[f_i(x)]^p$ ~where $x$ equals $n$ and $m$ respectively.
Because in (5,~5') each of the four $f_i()$ terms has form $bp^j+a$ mod
$p^{j+1}$ where the resp. $a<p^j$ yield (5), and the resp. msd's~ $b<p$
cause inequivalence (5'). ~~Then:\\[1ex] \hspace*{1cm}
(6)~~~~~ $f_{i+1}() \e (bp^j+a)^p \e a^{p-1}bp^{j+1}+a^p$ mod $p^{j+2} \e a^p$
        mod $p^{j+1}$ \\[1ex]
which depends only on $a$, and not on msd $bp^j$ of $f_i()$. This
preserves equivalence (5) mod $p^j$ for $i \rightarrow i+1$, and similarly
inequivalence (5') mod $p^{j+1}$ because, depending only on the respective
$a$ mod $p^j$, equivalence at $i+1$ would contradict (5') at $i$.
Cases $i<j$ and $i \geq j$ behave as follows. \\[1ex]
For $i<j$ the successive differences~
 $e_i(n)-e_i(m) \e y_ip^j \nqv 0$ mod $p^{j+1}$ . . . (6')  vary with $i$
 from 1 to $j-1$, and by (3') the core residues $f_i()$ mod $p^i$ settle for
increasing precision $i$. \\So initial inequivalences mod $p^p$ (4), and more
specifically mod $p^{j+1}$ (5), are preserved. \\And for all $i \geq j$ the
differences (6') are some constant $cp^j \nqv 0$ mod $p^{j+1}$, again by (3').
Hence by induction base (4) and steps (5,~6): core $A_p$ mod $p^p$ has
$h=(p-1)/2$ distinct increments, so critical precision $K_p<p$.
\end{proof}

Apparently $K_p$ is determined already by the initial integer increments
$e_1(n) < p^p ~(0<n<p)$, as the minimum precision $k$ for which non-symmetric
$n,m<p$ (so $n+m \neq p-1$) have $e_1(n) \nqv e_1(m)$ mod $p^k$.\\
For instance $p$=11 has $K_p=3$, and mod $p^3$ we have $h=5$ distinct core
increments, in base 11 code: $d_3(1..9)=\{4a1,~711,~871,~661,~061,~661,~871,
~711,~4a1\}$
~so core $A_3$ has the maximal five cosets generated by increments $d_3(n)$.
Equivalence $d_2(4) \e d_2(5) \e 61$ mod $p^2$ implies 661 and 061 to be in the
same $F$-coset in $G_3$. In fact 061.601=661 (base 11) with 601 in $F$
mod $p^3$, as are all $p$ residues of form $\{mp^2+1\} \e (p^2+1)^*$ mod $p^3$.

As example of lem3.1c, with $p=11$ and upto 3-digit precision:
\\ \hspace*{1in}
$\{n^p\}~~=~\{001,~5a2,~103,~274,~325,~886,~937,~aa8,~609,~0aa\}$\\
\hspace*{1in}
   core $A_3=\{001, 4a2,~103,~974,~525,~586,~137,~9a8,~609,~aaa\}$\\
$e_1(4)=325-274=061$ ~and \\$e_1(5)=886-325=561$
~with $FST$-carries: $4^{p-1}=a1,~5^{p-1}=71,~6^{p-1}=51$ ~so: \\[1ex]
    $e_2(4)=525-974=661$ ~by rule(3) yields:
    $5^{p^2}-4^{p^2}=[701]5^p-[a01]4^p=661$\\
    $e_2(5)=586-525=061$ derived by (3) as:
    $6^{p^2}-5^{p^2}=[501]6^p-[701]5^p=061$\\
Notice second difference $e_2(5)-e_2(4)=061-661=500$ equals
$e_1(5)-e_1(4)=561-061=500$ by lem3.1(c).

With $|F|=|G|/p$ and $|D_k|$ equal to $(p-1)/2$ for large enough $k<p$,
the nonzero $p$-th power pairsums cover nearly half of $G$. It will be shown
that four $p$-th power residues suffice to cover not only $G$ mod $p^k$, but
all residues $Z$ mod $p^k$. ~In this additive analysis we use:

{\bf Notation}: ~$S_{+t}$ is the set of all sums of $t$ elements in set $S$,
\\ \hspace*{1in} and $S+b$ stands for all sums $s+b$ with $s \in S$.

Extension subgroup $B$ is much less effective as additive generator than $F$.
Notice that $B \e \{np+1\}$ sothat $B+B \e \{mp+2\}$, and in general
$B_{+i} \e \{np+i\}$ in $G$, denoted by $N_i$, the subset of $G$ which is
$i$ mod $p$. They are also the (additive-) {\bf translations} $N_i \e
B-1+i ~(i<p)$ of $B$. Then $N_1 \e B$, while only $N_0 \e \{n.p\}$ is
not in $G$, and $N_i+N_j \e N_{i+j}$, corresponding to addition mod $p$.

Coresums $A_{+i}$ in general satisfy the next inclusions, implied by
$0 \in A_{+2} \e A+A$ :

For all $i \geq 1$: ~~$A_{+i} \subseteq A_{+(2+i)}$,
   ~and  $F_{+i} \subseteq F_{+(2+i)}$.

$F_{+3}$ covering all non-zero multiples $mp$ mod $p^k ~(k \geq 2)$ in
$N_0$ is related to a special result on the number 2 as generator.
For instance, a computer scan showed~ $2^p \nqv 2$ mod $p^2~(2 \notin A_2$)
for all primes $p<10^9$ except 1093 and 3511, although inequality does hold
mod $p^3$ for all primes (shown next). ~Notice that only 2 divides
$p-1$ for each odd prime $p$, so the 2-cycle $C_2 \e \pm 1$ is the only cycle
common to all cores for $p>$2. The generative power of 2 might be related
to it being a divisor of $p-1$ and $p+1$ for all $p>$2.

Regarding the known unsolved problem of a simple rule to find primitive
roots of 1 mod $p^k$, consider the divisors $r$ of $p^2-1=(p-1)(p+1)$ as
$generators$. \\
Recall that by (1) units group $G_k \e A_kB_k$ mod $p^k$ has core subgroup $A_k$
of order $p-1$, for any precision $k>0$, and extension group $B_k \e (p+1)^*$
of all $p^{k-1}$ residues 1 mod $p$, generated by $p+1$ ~[4,lem2]. In fact
$p-1$ generates all $2p^{k-1}$ residues $\pm 1$ mod $p^k$, including $B_k$.
\\
In multiplicative cyclic group $G_k$ of order $(p-1)p^{k-1}$, it stands to
reason to look for generators of $G_k$ (primitive roots of 1 mod $p^k$)
among the divisors of such powerful generators as $p \pm 1$, or similarly
of $p^2-1=(p-1)(p+1)$. Given prime structure $p^2-1=\prod_i ~p_i^{e_i}$
there are $\prod_i ~(e_i+1)$ divisors, forming a lattice, which is not
Boolean since factor $2^2$ makes $p^2-1$ non-squarefree.

Notice that for each unit $n$ in $G_k$ we have $n^{p-1}$ in $B_k$, and
$n^{p^{k-1}}$ in core $A_k$, while intersection $A_k ~\cap~ B_k  \e 1$ mod $p^k$,
the single unity of $G_k$. No generator $g$ of $G_k$ can be in core $A_k$,
since $|g^*|=(p-1)p^{k-1}$, while the order $|n^*|$ of $n \in A_k$ divides
$|A_k|=p-1$. ~Hence $p$ must divide the order of any non-core residue. If
$n<p^k$ then $n$ can be interpreted both as integer and as residue mod $p^k$.
It turns out that analysis modulo $p^3$ suffices to show that the divisors
$r$ of $p \pm 1$ are outside core, so $r^p \nqv r$ mod $p^3$: a necessary
but not sufficient condition for a primitive root. This amounts to quadratic
analysis of an extension of Fermat's Small Theorem ($FST$) on $p$-th power
residues, including two carry digits (base $p$).

\begin{thm} ~( divisors of $p^2-1$ ) \\
\hspace*{1in} If~ $r>1$ divides ~$p^2-1$
        ~~then ~$r^p \nqv r$ mod $p^k ~~~~(k \geq 3)$.
\end{thm}
\begin{proof}
~~~$r^p \nqv r$ mod $p^k$ implies inequality mod $p^{k+1}$. With $A_2 \e
F_2 \e \{n^p\}$ mod $p^2$, so each $p$-th power is in core $A_2$ mod $p^2$,
it suffices to show $r^p \nqv r$ ~mod $p^3$.
~Factorize $p^2-1=rs$, with positive integer cofactors $r$ and $s$. Then
$rs \e -1$ mod $p^2$, so opposite signed cofactors $\{r,-s\}$ or $\{-r,s\}$
form an inverse pair mod $p^2$. Inverses in a finite group $G$ have equal
order ($period$) in $G$, with order two automorphism $n \leftrightarrow
n^{-1}$. So orders $|r^*|$ and $|(-s)^*|$ are equal in $G_2$.
\\[1ex]
Notice $rs=p^2-1$ is not in core $A_3$, where $-1$ mod $p^3$ is the only
core residue that is $-1$ mod $p$, since the $p-1$ core residues $n^{|B_k|}$
of $A_k$ are distinct $\nqv 0$ mod $p ~(FST)$. In fact $(rs)^p \e (p^2-1)^p \e 1$
mod $p^3$ and no smaller exponent yields this. So $p^2-1=rs$ has order $2p$
in $G_3$, generating all $2p$ residues $\pm 1$ mod $p^2$, with inverse
pair $\{r^p,-s^p\}$ of equal order in $G_3$. Core $A_3$ is closed under
multiplication, so at most one cofactor of non-core product $rs$ can be
in core. In fact neither is in core $A_3$, so both $r^{p-1}$ and $s^{p-1}$
are $\nqv 1$ mod $p^3$, seen as follows.

 By $G_3=A_3B_3 ~(1)$: each $n \in G_3$ has product form $n \e n'.n"$
 mod $p^3$ of two components, with $n'$ in core $A_3$ and $n"$ in
 extension group $B_3$. Then $r^p(-s)^p \e 1$ mod $p^3$, where $r^p$ and
 $-s^p$ as inverse pair in $G_3$ have equal order, and each component
 forms an inverse pair of equal (and coprime) orders in $A_3$ and $B_3$
 respectively. The latter must divide $|B_3|=p^2$, and discarding order 1
(both $r,s$ cannot be in core, as shown) their common order is $p$ or $p^2$ ~[*].
For any unit $n$ the order of $n^p$ divides that of $n$, so $p$ dividing the
common order of $r^p$ and $s^p$ implies $p$ dividing also those of $r$ and $s$,
hence cofactors $r$ and $s$ of $p^2-1$ are both outside core $A_3$.
\end{proof}

PS [*]: "... and discarding order 1 ..." fails if $r^p \e r$ mod $p^2$,
which can occur (e.g. $p=11,~r=3$). This gap, noted by R. Chapman
(see his recent letter to the editor, and my reply), is mended by a proof
as given in a short paper "On primitive roots of unity,~ divisors of~
$p^2-1$, and an extension to mod $p^3$ ~of~ Fermat's Small Theorem"
- http://arXiv.org/abs/math.GM/0103067 submitted to {\it Computers and
Mathematics with Applications}. 

{\bf Notes}:
\begin{enumerate} {\small
\item A generator $g<p$ of $G_2$, so $|g^*|=(p-1)p$, also
    generates $G_k$ mod $p^{k>2}$ of order $(p-1)p^{k-1}$ ~[1].

\item Cofactors $r,s$ in $rs=p^2-1=(p-1)(p+1)$ have equal period in $G_3$,
    upto a factor of 2, so only $r \leq p+1$ need be inspected for
    periodic analysis. Recall exceptions $p=1093,~3511$ with $2^p=2$
    mod $p^2$, the only two primes $p<10^9$ with this property. Of the 79
    primes upto 401 there are seven primes with $r^p \e r$ mod $p^2$ for some
    divisor $r~|~p^2-1$ and cofactor $s$,  namely:\\ \hspace*{1cm}
    $p(r):~ 11(3),~29(14),~37(18),~181(78),~257(48),~281(20),~313(104)$.

\item A generator $g$ of $G_k$ is outside core, but $g~|~p \pm 1$ (thm3.1)
    does not guarantee $G_k \e g^*$.

\item However, computational evidence seems to suggest the next {\bf Conjecture}:\\
  At least one divisor $g~|~p \pm 1$ (prime $p>2$)
       generates $G_k$, or half of $G_k$ with $-1$ missing:\\
   then complements $-n$ mod $p^k$ yield the other half of $G_k$
	(e.g.~$p$=73:~ $G_3 \e \pm 6^* \e \pm 12^*)$.

\item Similarly the theorem also holds for divisors of $p^2+1$;
  and for $r~|~p^4-1=(p^2-1)(p^2+1)$, etc.: ~$r~|~p^{2m}-1$.
}
\end{enumerate}

For odd prime $p$ holds: ~2 divides both $p-1$ and $p$+1,
and 3 divides one of them, hence:

\begin{cor} 
For prime $p$ (incl. $p=2), ~k \geq 3$ and $n=2,~3$: \\
\hspace*{.5cm}
  $n^p \nqv n$ mod $p^k$, and in fact $\pm~\{n,~n^{-1}\}$ mod $p^k$
  are outside core $A_{k>2}$ for {\rm every odd prime}.
\end{cor}
In set notation: quadruple $Q(r)=\pm~\{r,~r^{-1}\},~r|(p \pm 1)$ and
~$k \geq 3$ ~imply ~$Q(r) \cap A_k=\emptyset$. \\
Moreopver, the product of $r \notin A_k$ with a core element is outside
core: ~$[~Q(r)A_k~] \cap A_k = \emptyset$. \\
Hence 2 is not in core $A$ mod $p^k$ for any prime $p>$2. This relates
to $p-1$ having divisor 2 for all $p$, and ~$C_2=\{-1,1\}$ as the only
common subgroup of $Z$(.) mod $p^k$ for all primes $p>$2. And 2 not in
core implies the same for its complement and inverse, $-2$ and $\pm~2^{-1}$.

Notice that $N_0$ mod $p^k$ consists of all multiples $mp$ of $p$, and their
base $p$ code ends on '0', so $|N_0|=p^{k-1}$. In fact $N_0$ consists of
all divisors of 0, the maximal nilpotent subsemigroup of $Z(.)$ mod $p^k$,
the semigroup of residue multiplication. For prime $p$ there are just
two idempotents in $Z(.)$ mod $p^k$: 1 in $G$ and 0 in $N_0$, so $G$ and
$N_0$ are complementary in $Z$, noted ~$N_0 \e  Z~\backslash~G$.

For prime $p>2$, consider integer $p$-th power function $F(n)=\{n^p\}$, with
$F_k$ denoting set $F(n)$ mod $p^k$ for all $n \nqv 0$ mod $p$, and core
function $A_k(n) \e n^{p^{k-1}}$, with core $A_2 \e F_2$. Multiples $mp$
($m \nqv 0$ mod $p$) are not $p$-th power residues (which are 0 mod $p^2$),
thus are not in $F_k$ for any $k>$1. But they are sums of three $p$-th power
residues: $mp \in F_{+3}$ mod $p^k$ for any $k>$1, shown next. In fact, due to
$FST$ we have $F(n) \e n$ mod $p$ for all $n$, so $F(r)+F(s)+F(t) \e r+s+t$ mod $p$,
which for a sum 0 mod $p$ of positive triple $r,s,t$ implies $r+s+t=p$.

\begin{lem} 
For $m \nqv 0$ mod $p$:~ $mp \in ~F_{+3}$ mod $p^{k>1}$, hence:\\
\hspace*{1cm} Each multiple $mp$ mod $p^{k>1}$ outside $F_k$
              is the sum of three $p$-th power residues (in $F_k$)
\end{lem}
\begin{proof}
Analysis mod $p^2$ suffices, because each $mp$ mod $p^{k>1}$ is reached upon
multiplication by $F_k$, due to (.) distributing over (+). Core $A_k$ has
order $p-1$ for any $k>$0, and $F_2 \e A_2$ implies powersums $F_2+F_2+F_2$ mod
$p^2$ to be sums of three core residues. \\
Assume $A(r)+A(s)+A(t) \e mp \nqv 0$ mod $p^2$ for some positive $r,s,t$ with
$r+s+t=p$. \\ Such $mp \notin A_2$ generates all $|A_2.mp|=|A_2|=p-1$ residues
in $N_0 ~\backslash ~0$ mod $p^2$. And for each prime $p>2$ there are many
such coresums $mp$ with $m \nqv 0$ mod $p$, seen as follows.

Any positive triple $(r,s,t)$ with $r+s+t=p$ yields, by $FST$, coresum
$A(r)+A(s)+A(t) \e r+s+t \e p$ mod $p$, hence with a coresum $mp$ mod $p^2$.
If $m$=0 then this solves $FLT~case_1$ for residues mod $p^2$, for instance
the cubic roots of 1 mod $p^2$ for each prime $p \e 1$ mod 6, see [4].\\
Non-zero $m$ is the dominant case for any prime $p>2$. In fact, normation
upon division by one of the three core terms in units group $G_2$ yields
one unity core term, say $A(t) \e 1$ mod $p^2$ hence $t=1$. Then $r+s=p-1$
yields $A(r)+A(s) \e mp-1$ mod $p^2$, where $0<m<p$. \\ There are
$1 \leq |D_2| \leq (p-1)/2$ distinct cosets of $F_2 \e A_2$ in $G_2$
(lem2.1, lem3.1), yielding as many distinct core pairsums $mp-1$
mod $p^2$ in set $A_2+A_2$.
\end{proof}

For most primes take $r=s$ equal to $h=(p-1)/2$ and $t=1$, with core residue
$A(h) \e h \e -2^{-1}$ mod $p$. Then $2.A(h)+1 \e mp \e 0$ mod $p$, with summation
indices $h+h+1=p$. For instance $p$=7 has $A(3) \e 43$ mod $7^2$ (base 7),
and $2A(3)+1=16+1=20$.

If for some prime $p$ we have in this case $m \e 0$ mod $p$, then $2.A(h) \e -1$
mod $p^2$, hence $A(h) \e h^p \e h$ mod $p^2$ and thus also $A(2) \e 2^p \e 2$ mod $p^2$.
In such rare cases (for primes $<10^9$ only $p=1093$ and $p=3511$) a choice of
other triples $r+s+t=p$ exists for which $A(r)+A(s)+A(t) \e mp \nqv 0$ mod $p^2$,
as just shown.\\
For instance $2^p \e 2$ mod $p^2$ for $p$=1093, but $3^p=936p+3$ mod $p^2$ sothat
instead of $(h,h,1)$ one applies $(r,s,1)$ where $r=(p-1)/3$ and $s=(p-1)2/3$.
And $p=3511$ has $3^p \e 21p+3$ mod $p^2$, while $3|p-1$ allows a similar index
triple with coresum $mp \nqv 0$ mod $p^2$.\\
Lemma 3.2 leads to the main additive result for residues in
ring $Z$[+, .] mod $p^k$ : \\[1ex] \hspace*{1cm}
   Each residue mod $p^k$ is the sum of at most four $p$-th power residues.

In fact, with ~subgroup $F \e \{n^p\}$ of $G$ in semigroup $Z$(.) mod $p^k$,
subsemigroup $N_0 \e \{mp\}$ of divisors of zero, and extension group
$B \e N_1 \e N_0$+1 in $G$, we have :

\begin{thm} 
~For residues mod $p^k ~(k \geq 2$, ~prime $p>$2):~~~
    $Z \e N_0~\cup~G \e F_{+3} ~\cup ~F_{+4}$
\end{thm}
\begin{proof}
~Analysis mod $p^2$ suffices, by extension lem3.1, and by lem3.2 all non-zero
multiples of $p$ are $N_0 ~\backslash~0 \e ~F_{+3}$, while $0 \in F_{+2}$
because $-1 \in F$. Hence $F_{+2} ~\cup ~F_{+3}$ covers $N_0$. Adding an extra
term $F$ yields $F_{+3} ~\cup ~F_{+4} ~\supseteq N_0+F$, which also covers
$A.N_0+A \supseteq A.(N_0+1) \e A.B \e G$ because $1 \in A$ and $A \subseteq F$,
so all of $Z \e N_0 \cup G$ is covered.
\end{proof}

{\bf Notes}:~ \begin{enumerate}  {\small
\item
 Case $p=3$ is easily verified by complete inspection as follows. Analysis
mod $p^3$ (thm.2) is rarely needed; for instance condition $2^p \nqv 2$ mod
$p^2$ holds for all primes $p<10^9$ except for the two primes 1093 and 3511.
So mod $p^2$ will suffice for $p$=3; moreover $F  \e  A$ mod $p^2$.  \\
Now $F \e \{-1,1\} \e \pm 1$ sothat $F+F \e \{0, \pm 2\}$.
Adding $\pm 1$ yields $F_{+3} \e \pm \{1,3\}$, and again $F_{+4} \e
\{0,\pm 2, \pm4\}$, sothat $F_{+3} \cup F_{+4}$ indeed cover all residues
mod $3^2$. Notice that $F_{+3}$ and $F_{+4}$ are disjoint which, although an
exception, necessitates their union in the general statement of thm3.2. \\
--- It is conjectured that $F_{+3} \subseteq F_{+4}$ for $p>$6,
    then $Z \e F_{+4}$ for primes $p>$6.

\item
For $p$=5 again use analysis mod $p^2$, and test if $F.(2A(h)+1)$ covers
all nonzero $m.5$ mod $5^2$ (lem3.2). Again $F \e A$ mod $p^2$, implying
$A(h) \in F$. Now core $A \e F \e (2^5)^* \e \{7, -1, -7, 1\}
\e \pm \{1, 7\}$, while $h \e 2$ with $A(2) \e 7$, or in base 5
code: $A(2) \e 12$ and $2A(h)+1 \e 30$. ~Hence $F(2A(h)+1) \e
\pm \{01, 12\}30 \e \pm \{30, 10\}$, indeed covering all
four nonzero residues $5m$ mod $5^2$.  }
\end{enumerate}

\section*{ Conclusions }

The application of elementary semigroup concepts to structure analysis of
residue arithmetic mod $p^k$ [3,4, 6] is very useful, allowing divisors of
zero. Fermat's inequality and Waring's representation are about powersums,
thus about additive properties of closures in $Z(.)$ mod $p^k$.

Fermat's inequality, viewed as anti-closure, reveals $n^p$ as a powerful
set of additive generators of $Z(+)$. Now $Z(.)$ has idempotent 1,
generating only itself, while 1 generates all of $Z(+)$ (Peano).

Similarly, expanding 1 to the subgroup $F \e \{ n^p \}$ of $p$-th power
residues in $Z(.)$ mod $p^k$, of order $|F|=|G|/p$, yields a most efficient
additive generator with: $F_{+3} \cup F_{+4} \e Z(+)$ mod $p^k$ for any
prime $p>2$. This is compatible for $p$=2 with the known result of each
positive integer being the sum of at most four squares.

The concept of $critical~precision$ (base $p$) is very useful for linking
integer symmetric properties to residue arithmetic mod $p^k$, and quadratic
analysis (mod $p^3$) for generative purposes such as primitive roots.
Eventually for binary arithmetic $p=2: ~p^2-1=p+1=3$ is a semi-
primitive root of 1 mod $2^k$ for $k \geq 3$ (thm3.1: note4, [4]: lem2)
with useful engineering result [8].

\section*{ References }
\begin{enumerate}
\item T.Apostol: {\it Introduction to Analytical Number Theory},
  thm 10.4-6, Springer Verlag 1976
\item E.T.Bell: {\it The development of mathematics}
   ~(p304-6) McGraw-Hill, 1945
\item N.F.Benschop: "The semigroup of multiplication mod $p^k$, an
  extension of Fermat's Small Theorem, and its additive structure"
  {\it Semigroups and Applications}, p7, Prague, July'96.
\item --- : "The triplet structure of arithmetic mod $p^k$, Fermat's
  Small and Last Theorem, and a new binary number code."
  {\it Logic and Architecture Synthesis},  p133-140, Grenoble, Dec'96.
  ~(also as http://www.iae.nl/users/benschop/199706-1.dvi )
\item A.Clifford, G.Preston: {\it The Algebraic Theory of Semigroups} \\
   Vol 1 (p130-135), AMS survey \#7, 1961.
\item S.Schwarz: "The Role of Semigroups in the Elementary Theory of
   Numbers",\\ Math.Slovaca V31, N4, p369-395, 1981.
\item G.Hardy, E.Wright: {\it An Introduction to the Theory of Numbers} \\
  (Chap 8.3, Thm 123), Oxford-Univ. Press 1979.
\item N.F.Benschop: ~Patent US-5923888 Logarithmic multiplier
  (dual bases 2 and 3) 13jul99.
\end{enumerate}

\end{document}